\documentclass[11pt]{amsart}
\usepackage[T1]{fontenc}
\usepackage[utf8]{inputenc}
\usepackage[english]{babel}
\usepackage{amsmath,amssymb,amsthm,mathtools}
\usepackage[pdftex,unicode]{hyperref}
\usepackage{enumitem}
\setlist[enumerate,1]{label=\textup{(\roman*)}}

\newtheorem{theorem}{Theorem}[section]
\newtheorem{lemma}[theorem]{Lemma}
\newtheorem{proposition}[theorem]{Proposition}
\newtheorem{cor}[theorem]{Corollary}
\theoremstyle{definition}

\newtheorem{example}[theorem]{Example}
\newtheorem{question}[theorem]{Question}
\newtheorem{problem}{Problem}
\theoremstyle{remark}
\newtheorem{note}[theorem]{Remark}

\def\R{\mathbb R}
\def\N{\mathbb N}
\def\I{\mathbb I}
\def\om{\omega}
\def\si{\sigma}
\def\de{\delta}
\def\ep{\varepsilon}

\def\la{\lambda}
\def\be{\beta}

\def\es{\varnothing}

\def\cA{\mathcal A}

\def\cK{\mathcal K}

\def\cl#1{\overline{#1}}
\def\clx#1#2{\overline{#2}^{\,#1}}

\def\set#1{\bbset#1\eeset}
\def\bbset#1:#2\eeset{\{#1\,:\,#2\}}

\def\:{\colon}
\def\Hc{\I^\om}
\def\fb{\bar f}
\def\fs{f_\sharp}

\def\pq{\preccurlyeq}
\def\ph{\varphi}

\begin{document}

\title[$C$-embeddings of spaces of signed measures]{$C$-embeddings of spaces of signed measures}

\author{E.~A.~Reznichenko}
\address{Department of Mathematics, Lomonosov Moscow State University, Moscow, Russia}
\author{Yu.~V.~Sadovnichiy}
\address{Department of Mathematics, Lomonosov Moscow State University, Moscow, Russia}

\begin{abstract}
For a Tychonoff space $X$ let $M_a(X)\supset M_\si(X)\supset M_\tau(X)\supset M_R(X)$ be the spaces of finitely additive, $\si$-additive (Baire), $\tau$-additive and Radon signed measures on $X$, and let $M_\be(X)$ and $M_d(X)$ be the spaces of compactly supported and discrete measures; all of them are regarded as subspaces of $M(\be X)=C(\be X)^*$ with the weak$^*$ topology. We study when $M_p(X)$ is $C$-embedded ($C^*$-embedded) in $M_q(X)$. We prove that $C$-embeddedness, $C^*$-embeddedness and $G_\de$-density of $M_p(X)$ in $M_q(X)$ are equivalent, and that they are equivalent to a lifting property: for every continuous map $f$ of $X$ into a metrizable compact space, the image of every probability measure in $M_q(X)$ is the image of a probability measure in $M_p(X)$. The main tool is a general theorem, proved in a separate paper, asserting that $C$-embeddedness, $C^*$-embeddedness and $G_\de$-density are equivalent for dense subspaces of locally convex spaces in the weak topology, and, more generally, of products of first-countable spaces. We characterize the spaces $X$ for which $M_\si(X)$, $M_\be(X)$, $M_d(X)$ and $M_R(X)$ are $C$-embedded in $M_a(X)$ (pseudocompactness, pseudo-$\om$-boundedness, scatteredness of $\be X$, and a lifting property of Radon measures, respectively). We give an example of a pseudocompact locally compact space $X$ for which $M_\tau(X)$ is not $C^*$-embedded in $M_\si(X)$, and we show that $M_\tau(X)$ is $C$-embedded in $M(\be X)$ for every $G_\de$-dense subspace $X$ of a product of metrizable compact spaces. Using a set of universal measure zero, we construct a pseudocompact space $X$ for which $M_R(X)$ is $C$-embedded in $M(\be X)$, although the restriction to $X$ of a coordinate projection onto $[0,1]$ maps no $\si$-compact subset of $X$ onto $[0,1]$.
\end{abstract}


\maketitle

%======================================================================
\section{Introduction}\label{sec:intro}

All spaces are assumed to be Tychonoff. A subspace $Y$ of a space $Z$ is $C$-embedded ($C^*$-embedded) in $Z$ if every continuous (bounded continuous) function $Y\to\R$ extends to a continuous function on $Z$.

Let $X$ be a space. As in \cite{RS}, we identify the classical spaces of signed measures on $X$ with subspaces of $M(\be X)$, the dual of the Banach space $C(\be X)$ with the weak$^*$ topology:
\begin{gather*}
M_d(X)\subset M_R(X),\\
M_\be(X)\subset M_R(X)\subset M_\tau(X)\subset M_\si(X)\subset M_a(X)=M(\be X)
\end{gather*}
(the precise definitions are recalled in Section~\ref{sec:prel}). We denote by $M^+_p(X)$ and $P_p(X)$ the corresponding subsets of nonnegative and probability measures; for a continuous map $f$ of $X$ into a metrizable compact space $K$, $f_\sharp\:M(\be X)\to M(K)$ is the induced map (Section~\ref{sec:prel}).

This work was motivated by several questions on spaces of probability measures.

Banakh, Chigogidze and Fedorchuk \cite{BanakhChigogidzeFedorchuk2003} proved that the space $P_\si(X)$ is realcompact. The spaces $P_R(X)\subset P_\tau(X)$ are dense in $P_\si(X)$. The authors of \cite{BanakhChigogidzeFedorchuk2003} studied when $P_\si(X)$ is the Hewitt realcompactification of $P_\tau(X)$ or of $P_R(X)$. In other words, they asked when the functors $P_\si$ and $\nu\circ P_\tau$ (respectively, $\nu\circ P_R$) are isomorphic. They found conditions under which the answer is positive \cite[Corollaries 5.2 and 5.3]{BanakhChigogidzeFedorchuk2003}. A dense subspace $Y$ of a realcompact space $Z$ satisfies $\nu Y=Z$ if and only if $Y$ is $C$-embedded in $Z$. Thus the question is a question on $C$-embeddings.

\begin{enumerate}[leftmargin=3.5em]
\item[(BCF)]
Is $P_\si(X)$ the Hewitt realcompactification of $P_\tau(X)$ for every space $X$? Equivalently, is $P_\tau(X)$ $C$-embedded in $P_\si(X)$? \cite[Problem 5.5]{BanakhChigogidzeFedorchuk2003}
\end{enumerate}

Bogachev \cite{Bogachev2024-1ru,Bogachev2024dan-ru} studied the commutation of the functors $\be$ and $P=P_R$, that is, the equality $\be P(X)=P(\be X)$. This equality means that $P(X)$ is $C^*$-embedded in $P(\be X)$. Bogachev proved that the equality implies the pseudocompactness of both $X$ and $P(X)$. He also found spaces for which the equality holds; examples are $\om_1$ and the Tychonoff plank. He asked whether the converse implications hold.

\begin{enumerate}[leftmargin=3.5em]
\item[(B1)]
Does pseudocompactness of $X$ imply $\be P(X)=P(\be X)$? \cite{Bogachev2024-1ru}
\item[(B2)]
Does pseudocompactness of $P(X)$ imply $\be P(X)=P(\be X)$? \cite{Bogachev2024-1ru}
\end{enumerate}

These questions were answered by the first author in \cite{Reznichenko2025betap}. The answer to (B2) is positive: $P(X)$ is $C^*$-embedded in $P(\be X)$ if and only if $P(X)$ is pseudocompact \cite[Theorem 1.9]{Reznichenko2025betap}. The answers to (BCF) and (B1) are negative. There exist pseudocompact spaces $X$ for which $P(X)$ is not pseudocompact \cite[Examples 1.13 and 1.15]{Reznichenko2025betap}. One of them is a Mr\'owka--Isbell space $X$ with $|\be X\setminus X|=1$; it is pseudocompact and locally compact. Following \cite{Reznichenko2025betap}, we call this particular space the Mr\'owka space, and use the term Mr\'owka--Isbell space for the general class. For this space $P_\tau(X)=P(X)$ is not $C$-embedded in $P_\si(X)=P(\be X)$ \cite[Example 1.23]{Reznichenko2025betap}.

The analogous questions for signed measures are the following. Is $M_\tau(X)$ $C$-embedded in $M_\si(X)$? Does pseudocompactness of $X$, or of $P(X)$, imply that $M_R(X)$ is $C$-embedded in $M(\be X)$? We show that the answers do not reduce to the case of probability measures. As in the probability case, pseudocompactness of $X$ is not sufficient (Example~\ref{ex:psi}), while pseudocompactness of $P(X)$ is sufficient (Theorem~\ref{t:R}). However, for the Mr\'owka space of \cite{Reznichenko2025betap} the space $M_R(X)$ is $C$-embedded in $M(\be X)$, and $M_\tau(X)$ is $C$-embedded in $M_\si(X)$. At the same time, $P(X)$ is not pseudocompact, and $P_\tau(X)$ is not $C$-embedded in $P_\si(X)$ (Example~\ref{ex:mrowka}). Consequently, Bogachev's theorem does not extend to signed measures in full. If $M_R(X)$ is $C$-embedded in $M(\be X)$, then $X$ is pseudocompact (Corollary~\ref{c:pc}), but $P(X)$ need not be pseudocompact.

In \cite{RS} these questions were put into a general framework \textup(see also \cite[Problems 1.21 and 1.29]{Reznichenko2025betap}\textup). Let $p,q\in\{d,\be,R,\tau,\si,a\}$ and $Q\in\{M,M^+,U,P\}$, where $U_p(X)$ denotes the set of measures in $M_p(X)$ of norm at most one. A space $X$ is called universally $Q_{[p,q]}$-measurable (universally $Q^*_{[p,q]}$-measurable) if $Q_p(X)$ is $C$-embedded ($C^*$-embedded) in $Q_q(X)$. This is a functional version of the classical notions of universally measurable sets and sets of universal measure zero. For $X\subset\R$ these classical notions mean $M_\si(X)=M_R(X)$ and $M_\si(X)=M_d(X)$, respectively. We refer to \cite{RS} for this approach and use only its terminology. In the present paper we concentrate on the signed case $Q=M$; the cases $Q\in\{M^+,U,P\}$ are treated in the companion paper \cite{RS2}.

\begin{problem}\label{pr:main}
For which spaces $X$ is $M_p(X)$ $C$-embedded \textup($C^*$-embedded\textup) in $M_q(X)$?
\end{problem}

The basic tool is the following theorem on dense subsets of products of first-countable spaces, proved in the separate paper \cite{RS0}.

\begin{theorem}[{\cite[Theorem 1]{RS0}}]\label{t:intro:prod}
Let $\Pi$ be a product of first-countable spaces, $Z$ a dense subspace of $\Pi$, and $A$ a dense subset of $Z$. Then $A$ is $C$-embedded in $Z$ iff $A$ is $C^*$-embedded in $Z$ iff $A$ is $G_\de$-dense in $Z$ iff $\pi_J(A)=\pi_J(Z)$ for every countable set $J$ of coordinates, where $\pi_J$ is the projection onto the corresponding subproduct.
\end{theorem}

Every locally convex space in the weak topology is a dense linear subspace of a power of the real line. Therefore the same equivalence holds for dense subsets of such spaces (\cite[Theorem 2]{RS0}, see Theorem~\ref{t:lcs:main}). Applied to $M_q(X)\subset M(\be X)$, this gives the first part of the following general criterion (Theorem~\ref{t:main} and Proposition~\ref{p:MP}).

\begin{theorem}\label{t:intro:main}
Let $X$ be a space and $p\pq q$. The following conditions are equivalent:
\begin{enumerate}
\item
$M_p(X)$ is $C$-embedded in $M_q(X)$;
\item
$M_p(X)$ is $C^*$-embedded in $M_q(X)$;
\item
$M_p(X)$ is $G_\de$-dense in $M_q(X)$;
\item
$P_p(X)$ is $G_\de$-dense in $P_q(X)$;
\item
for every metrizable compact space $K$ and every continuous map $f\:X\to K$, every image $\fs\mu$ of a measure $\mu\in P_q(X)$ is the image $\fs\nu$ of a measure $\nu\in P_p(X)$.
\end{enumerate}
\end{theorem}

Thus, as far as $G_\de$-density is concerned, there is no difference between signed and probability measures. For $C^*$-embeddings there is a difference: for the Mr\'owka space of Example~\ref{ex:mrowka}, the space $P_R(X)$ is $G_\de$-dense, but not $C^*$-embedded, in $P(\be X)$. The spaces of nonnegative measures, probability measures and unit balls are studied in \cite{RS2}.

For $p\in\{d,\be,R\}$ condition (v) can be expressed without lifting of measures (Theorem~\ref{t:dbR}). As a consequence, for subsets of the real line the functional notions coincide with the classical ones (Corollary~\ref{c:classical}).

\begin{theorem}\label{t:intro:classical}
Let $X$ be a separable metrizable space. Then $M_R(X)$ is $C$-embedded in $M_\si(X)$ iff $X$ is universally measurable \textup(i.e.\ $M_R(X)=M_\si(X)$\textup), and $M_d(X)$ is $C$-embedded in $M_\si(X)$ iff $X$ has universal measure zero \textup(i.e.\ $M_d(X)=M_\si(X)$\textup).
\end{theorem}

For embeddings into $M_a(X)=M(\be X)$ we obtain the following characterizations (Corollary~\ref{c:pc} and Theorems~\ref{t:beta}, \ref{t:d}, \ref{t:R}).

\begin{theorem}\label{t:intro:a}
Let $X$ be a space.
\begin{enumerate}
\item
$M_\si(X)$ is $C$-embedded in $M(\be X)$ iff $X$ is pseudocompact.
\item
$M_\be(X)$ is $C$-embedded in $M(\be X)$ iff $X$ is pseudo-$\om$-bounded.
\item
$M_d(X)$ is $C$-embedded in $M(\be X)$ iff $\be X$ is scattered.
\item
$M_R(X)$ is $C$-embedded in $M(\be X)$ iff $X$ is pseudocompact and every Radon probability measure on a metrizable compact continuous image $f(X)$ of $X$ is the image of a Radon probability measure on $X$. This is the case if $P_R(X)$ is pseudocompact.
\end{enumerate}
\end{theorem}

Moreover, $M_\si(X)$ is $G_\de$-closed in $M(\be X)$, so for $p\pq\si$ the space $M_p(X)$ is $C$-embedded in $M(\be X)$ iff $X$ is pseudocompact and $M_p(X)$ is $C$-embedded in $M_\si(X)$. In (iv), the lifting condition cannot be replaced by the requirement that every continuous image $f(X)\subset\Hc$ be the image $f(S)$ of a $\si$-compact set $S\subset X$: this requirement is sufficient (Proposition~\ref{p:S}), but not necessary (Example~\ref{ex:83}).

Section~\ref{sec:example} contains examples. Let $X=\Psi(\cA)$ be the Mr\'owka--Isbell space of a suitable maximal almost disjoint family $\cA$; here $\be X$ is not scattered. Then $M_\tau(X)$ is not $C^*$-embedded in $M_\si(X)=M(\be X)$ (Example~\ref{ex:psi}). Thus the answers to the questions of \cite{BanakhChigogidzeFedorchuk2003} and \cite{Bogachev2024-1ru} remain negative for signed measures. For the $\Sigma$-product $X\subset[0,1]^{\om_1}$ all the inclusions $M_R(X)\subset M_\tau(X)\subset M_\si(X)$ are proper, and each of these spaces is $C$-embedded in each larger one (Example~\ref{ex:sigma}). The Mr\'owka--Isbell space $X$ of the Cantor tree is locally compact, submetrizable and realcompact, but $M_\tau(X)$ is not $C^*$-embedded in $M_\si(X)$ (Example~\ref{ex:tree}). On the other hand, for the realcompact space $X=\R^{\mathfrak c}$ we have $M_\tau(X)\ne M_\si(X)$, while $M_R(X)$ and $M_\tau(X)$ are $C$-embedded in $M_\si(X)$ (Remark~\ref{r:Rc}). If $X$ is a $G_\de$-dense subspace of a product of metrizable compact spaces, then $M_\tau(X)$ is $C$-embedded in $M(\be X)$ (Proposition~\ref{p:tau}). Let $X$ be a pseudocompact dense subspace of the Tychonoff or Cantor cube of weight $\mathfrak c$ all of whose compact subsets are countable; such spaces exist in ZFC, for example, the spaces of \cite{Shakhmatov1986,Reznichenko1989} and pseudocompact Korovin orbits \cite{Korovin1992,ReznichenkoTkachenko2024}. Then $M_\tau(X)$ is $C$-embedded in $M(\be X)$, but $M_R(X)$ is not $C$-embedded even in $M_\tau(X)$ (Example~\ref{ex:cube}). Finally, using an uncountable set of universal measure zero, we construct a pseudocompact space $X$ for which $M_R(X)$ is $C$-embedded in $M(\be X)$, while the restriction to $X$ of the projection onto the first factor $[0,1]$ maps no $\si$-compact subset of $X$ onto $[0,1]$ (Example~\ref{ex:83}).

The present paper continues \cite{RS}. In \cite{RS} the approach was outlined, and the main results were announced with sketches of proofs. Here we give detailed proofs. Theorem~\ref{t:main} below is a refined version of the main theorem of \cite{RS}. Its condition (v) is stated in terms of images of measures, and in this form the equivalence holds for all spaces (Remark~\ref{r:RS}). Accordingly, Corollaries~1 and~2 of \cite{RS} require additional assumptions on $X$. They hold, for example, for pseudo-$\om$-bounded spaces and for the Mr\'owka space of \cite{Reznichenko2025betap}. They do not hold for all pseudocompact spaces (Examples~\ref{ex:psi} and~\ref{ex:mrowka}).

%======================================================================
\section{Preliminaries}\label{sec:prel}

\subsection{Notation}
$\om=\{0,1,2,\dots\}$; $\Hc=[0,1]^\om$ is the Hilbert cube. A space is a metrizable compact space iff it is homeomorphic to a closed subspace of $\Hc$ \cite[Theorem 4.2.10]{EngelkingBookRu}; in proofs we shall often assume that a metrizable compact space is a subspace of $\Hc$. For a space $X$, $\be X$ is its Stone--\v Cech compactification; every bounded continuous function $h$ on $X$ is identified with its extension $\be h\in C(\be X)$. A zero set in a space $Z$ is a set of the form $h^{-1}(0)$, $h\in C(Z)$. In a compact space the zero sets are exactly the closed $G_\de$-sets.

A subset $S$ of a space $Z$ is $G_\de$-dense in $Z$ if $S$ meets every nonempty $G_\de$-subset of $Z$. If $S\subset T\subset Z$, then $S$ is $G_\de$-dense in $T$ iff $S$ meets every $G_\de$-subset of $Z$ that meets $T$.

\begin{lemma}\label{l:pc}
A space $X$ is pseudocompact iff every nonempty zero set of $\be X$ meets $X$.
\end{lemma}
\begin{proof}
If $Z=h^{-1}(0)\ne\es$, $h\in C(\be X)$, $h\ge 0$, and $Z\cap X=\es$, then $1/h$ is an unbounded continuous function on $X$ (since $X$ is dense in $\be X$). Conversely, if $g\:X\to[0,\infty)$ is continuous and unbounded, then $h=\be(1/(1+g))$ vanishes at some point of $\be X$ and does not vanish on $X$.
\end{proof}

A space $X$ is pseudo-$\om$-bounded if for every sequence $(U_n)_n$ of nonempty open subsets of $X$ there is a compact $K\subset X$ with $K\cap U_n\ne\es$ for all $n$ \cite{AngoaOrtiz-CastilloTamariz-Mascarua2014}. A space is scattered if each of its nonempty subsets has an isolated point. We shall use the following classical facts.

\begin{enumerate}
\item[(PS)]
A compact space $K$ is not scattered iff there is a continuous map of $K$ onto $[0,1]$ \cite{PelczynskiSemadeni1959}.
\item[(RK)]
Every Radon measure on a scattered compact space is discrete (Rudin \cite{Rudin1957}; see also \cite{Knowles1967}). An uncountable metrizable compact space contains a Cantor set and hence carries a nonzero Radon measure vanishing on points.
\end{enumerate}

\begin{lemma}\label{l:scat}
If $\be X$ is scattered, then $X$ is pseudocompact.
\end{lemma}
\begin{proof}
Let $g\:X\to\R$ be continuous and unbounded. Choose $x_k\in X$ with $|g(x_{k+1})|>|g(x_k)|+2$, and put $t_k=|g(x_k)|$, $D=\set{x_k:k\in\om}$. For a bounded sequence $(a_k)$ let $h\:[0,\infty)\to\R$ be the piecewise linear function with $h(t_k)=a_k$ that is constant on $[0,t_0]$. Then $h\circ|g|$ is a bounded continuous function on $X$ equal to $a_k$ at $x_k$. Hence the discrete space $D$ is $C^*$-embedded in $X$, and so in $\be X$. Therefore $\clx{\be X}D$ is homeomorphic to $\be\N$, which is not scattered, since $\be\N\setminus\N$ has no isolated points.
\end{proof}

\subsection{Measures on compact spaces}
For a compact space $K$ let $M(K)=C(K)^*$ with the weak$^*$ topology. By the Riesz theorem, $M(K)$ is the space of (signed) Radon measures on $K$; we write $\mu(\ph)=\int\ph\,d\mu$ for bounded Borel $\ph$. For $\mu\in M(K)$, $|\mu|$ is its total variation, $\mu=\mu^+-\mu^-$ its Jordan decomposition and $\|\mu\|=|\mu|(K)$. Let $M^+(K)$ be the cone of nonnegative measures and $P(K)=\set{\mu\in M^+(K):\mu(K)=1}$. We say that $\mu$ is concentrated on a Borel set $A\subset K$ if $|\mu|(K\setminus A)=0$.

For a continuous map $g\:K\to L$ of compact spaces, $g_\sharp=M(g)\:M(K)\to M(L)$, $g_\sharp\mu(\ph)=\mu(\ph\circ g)$, is linear and continuous. If $L$ is metrizable, the formula $g_\sharp\mu(\ph)=\mu(\ph\circ g)$ holds for all bounded Borel $\ph$ on $L$, because $B\mapsto\mu(g^{-1}(B))$ is a finite Borel (hence Radon) measure on $L$ that agrees with $g_\sharp\mu$ on $C(L)$. Clearly $g_\sharp$ maps $M^+(K)$ into $M^+(L)$ and $P(K)$ into $P(L)$, and $g_\sharp\mu$ is concentrated on $g(K)$.

\begin{lemma}\label{l:onto}
Let $K$ be a compact space and $D\subset K$ a dense subset.
\begin{enumerate}
\item
For $Q\in\{M,M^+,P\}$, the set of finitely supported measures $\sum_{i\le n} a_i\de_{x_i}\in Q(K)$ with $x_i\in D$ is dense in $Q(K)$.
\item
If $g\:K\to L$ is a continuous surjection of compact spaces, then $g_\sharp(Q(K))=Q(L)$ for $Q\in\{M,M^+,P\}$.
\end{enumerate}
\end{lemma}
\begin{proof}
(i) Let $C$ be the closure of the convex hull of $\set{\de_x:x\in D}$. Since $x\mapsto\de_x$ is continuous, $\de_y\in C$ for all $y\in K$. If $\mu\in P(K)\setminus C$, then by the separation theorem in the locally convex space $M(K)$ there are $h\in C(K)$ and $c\in\R$ with $\mu(h)>c\ge\nu(h)$ for all $\nu\in C$; taking $\nu=\de_y$ we get $h\le c$, hence $\mu(h)\le c$, a contradiction. So $C\supset P(K)$, which gives (i) for $Q=P$. The other cases follow, because the maps $(t,\mu)\mapsto t\mu$ and $(s,t,\mu,\nu)\mapsto s\mu-t\nu$ are continuous and $M^+(K)=\set{t\mu:t\ge 0,\ \mu\in P(K)}$ and $M(K)=\set{s\mu-t\nu:\ s,t\ge0,\ \mu,\nu\in P(K)}$.

(ii) The set $g_\sharp(P(K))$ is compact and contains all finitely supported probability measures on $L$, so it equals $P(L)$ by (i). The other cases follow from the representations of $M^+$ and $M$ through $P$ given above.
\end{proof}

\subsection{Spaces of measures on a Tychonoff space}
Let $X$ be a space. Following \cite{RS} (see also \cite{Wheeler1983,Bogachev2007Measure_theory}), for $p\in\{d,\be,R,\tau,\si,a\}$ the space $M_p(X)$ is the set of $\mu\in M(\be X)$ satisfying the condition $(M_p)$:
\begin{enumerate}
\item[$(M_d)$] $\mu$ is concentrated on a countable subset of $X$;
\item[$(M_\be)$] $\mu$ is concentrated on a compact subset of $X$;
\item[$(M_R)$] $\mu$ is concentrated on a $\si$-compact subset of $X$;
\item[$(M_\tau)$] $|\mu|(K)=0$ for every compact $K\subset\be X\setminus X$;
\item[$(M_\si)$] $|\mu|(Z)=0$ for every zero set $Z$ of $\be X$ with $Z\subset \be X\setminus X$;
\item[$(M_a)$] no condition, i.e.\ $M_a(X)=M(\be X)$.
\end{enumerate}
Measures in $M_d(X)$ are called discrete; thus a discrete measure is a measure concentrated on an at most countable subset of $X$ (not a measure on a discrete space). Under the standard identifications, $M_\si(X)$, $M_\tau(X)$, $M_R(X)$ are the spaces of countably additive Baire measures, $\tau$-additive Borel measures, and Radon measures on $X$, and $M_a(X)$ corresponds to finitely additive measures (see \cite{RS}); we shall not use these identifications.

We write $p\pq q$ if $M_p(X)\subset M_q(X)$ for all $X$; thus $d\pq R$, $\be\pq R\pq\tau\pq\si\pq a$, and $d$, $\be$ are incomparable. For $Q\in\{M^+,P\}$ let $Q_p(X)=M_p(X)\cap Q(\be X)$.

\begin{lemma}\label{l:Mp}
Let $p\in\{d,\be,R,\tau,\si,a\}$.
\begin{enumerate}
\item
$M_p(X)$ is a linear subspace of $M(\be X)$. If $\mu\in M_p(X)$, $\nu\in M(\be X)$ and $|\nu|\le|\mu|$, then $\nu\in M_p(X)$. In particular, $\mu^+,\mu^-,|\mu|\in M_p(X)$ for $\mu\in M_p(X)$.
\item
$\de_x\in M_d(X)\cap M_\be(X)$ for $x\in X$, and $Q_p(X)$ is dense in $Q(\be X)$ for $Q\in\{M,M^+,P\}$.
\end{enumerate}
\end{lemma}
\begin{proof}
(i) follows from $|\mu+\nu|\le|\mu|+|\nu|$, the fact that the conditions $(M_p)$ are conditions on $|\mu|$ that are inherited by smaller nonnegative measures, and from the fact that finite unions of countable (compact, $\si$-compact) sets are countable (compact, $\si$-compact). (ii) follows from Lemma~\ref{l:onto}(i) with $D=X$.
\end{proof}

For a metrizable compact space $K$ and a continuous map $f\:X\to K$ we denote by $\fb\:\be X\to K$ its continuous extension and put $\fs=\fb_\sharp\:M(\be X)\to M(K)$. For a bounded Borel function $\ph$ on $K$ we have $\fs\mu(\ph)=\mu(\ph\circ\fb)$. Note that $\fb(\be X)=\cl{f(X)}$, and $\fb(\be X)=f(X)$ if $X$ is pseudocompact.

\subsection{Dense subsets of products and of locally convex spaces}

We use two results of \cite{RS0}. For a set $S$ and spaces $X_s$, $s\in S$, let $\Pi=\prod_{s\in S}X_s$, and for $J\subset S$ let $\pi_J\:\Pi\to\Pi_J=\prod_{s\in J}X_s$ be the projection.

\begin{theorem}[{\cite[Theorem 1]{RS0}}]\label{t:prod:main}
Let $\Pi=\prod_{s\in S}X_s$ be a product of first-countable spaces, let $Z$ be a dense subspace of $\Pi$, and let $A$ be a dense subset of $Z$. The following are equivalent:
\begin{enumerate}
\item
$A$ is $C^*$-embedded in $Z$;
\item
$A$ is $C$-embedded in $Z$;
\item
$A$ is $G_\de$-dense in $Z$;
\item
$\pi_J(A)=\pi_J(Z)$ for every countable $J\subset S$.
\end{enumerate}
\end{theorem}

Let $L$ be a real vector space and $F$ a linear space of linear functionals on $L$ separating the points of $L$. Denote by $\si(L,F)$ the weakest topology on $L$ in which all $g\in F$ are continuous. The spaces $(L,\si(L,F))$ are exactly the locally convex spaces in the weak topology. The main example for us is $L=M(\be X)$, $F=\set{\mu\mapsto\mu(h):h\in C(\be X)}$: then $\si(L,F)$ is the weak$^*$ topology.

\begin{theorem}[{\cite[Theorem 2]{RS0}}]\label{t:lcs:main}
Let $Z$ be a dense subspace of $(L,\si(L,F))$ and $A$ a dense subset of $Z$. Then $A$ is $C^*$-embedded in $Z$ iff $A$ is $C$-embedded in $Z$ iff $A$ is $G_\de$-dense in $Z$.
\end{theorem}

Theorem~\ref{t:lcs:main} is deduced from Theorem~\ref{t:prod:main}, since $(L,\si(L,F))$ is linearly homeomorphic to a dense linear subspace of $\R^\Gamma$, where $\Gamma$ is a Hamel basis of $F$ \cite[Lemma 1]{RS0}.

%======================================================================
\section{The general criterion}\label{sec:main}

\begin{theorem}\label{t:main}
Let $X$ be a space and $p\pq q$. The following conditions are equivalent:
\begin{enumerate}
\item[(i)]
$M_p(X)$ is $C$-embedded in $M_q(X)$, i.e.\ $X$ is universally $M_{[p,q]}$-measurable;
\item[(ii)]
$M_p(X)$ is $C^*$-embedded in $M_q(X)$;
\item[(iii)]
$M_p(X)$ is $G_\de$-dense in $M_q(X)$;
\item[(iv)]
for every $\mu\in M_q(X)$ and every sequence $(h_n)_n\subset C(\be X)$ there is $\nu\in M_p(X)$ with $\nu(h_n)=\mu(h_n)$ for all $n$;
\item[(v)]
for every metrizable compact space $K$ and every continuous map $f\:X\to K$,
\[
\fs(M_q(X))=\fs(M_p(X)).
\]
\end{enumerate}
\end{theorem}
\begin{proof}
(i)$\Leftrightarrow$(ii)$\Leftrightarrow$(iii). Apply Theorem~\ref{t:lcs:main} to $L=M(\be X)$ with the weak$^*$ topology $\si(L,C(\be X))$, $A=M_p(X)$ and $Z=M_q(X)$; by Lemma~\ref{l:Mp}(ii), $A$ is dense in $L$, so $Z$ is dense in $L$ and $A$ is dense in $Z$.

(iii)$\Rightarrow$(iv). The set $\set{\nu\in M_q(X):\nu(h_n)=\mu(h_n)\text{ for all }n}$ is a $G_\de$-subset of $M_q(X)$ containing $\mu$.

(iv)$\Rightarrow$(iii). Let $H$ be a nonempty $G_\de$-subset of $M_q(X)$ and $\mu\in H$. Then $H=\bigcap_n O_n$, where $O_n$ is open in $M_q(X)$, and each $O_n$ contains a set $\set{\nu\in M_q(X):|\nu(h)-\mu(h)|<\ep_n,\ h\in A_n}$ with $A_n\subset C(\be X)$ finite. Enumerate $\bigcup_nA_n=\{h_k\}_k$; by (iv) there is $\nu\in M_p(X)$ with $\nu(h_k)=\mu(h_k)$ for all $k$, and then $\nu\in H$.

(iv)$\Rightarrow$(v). Let $f\:X\to K$ be continuous, $\mu\in M_q(X)$, and let $(g_n)_n$ be a dense sequence in the separable Banach space $C(K)$. By (iv) there is $\nu\in M_p(X)$ with $\nu(g_n\circ\fb)=\mu(g_n\circ\fb)$, i.e.\ $\fs\nu(g_n)=\fs\mu(g_n)$, for all $n$. Two elements of $M(K)$ coinciding on a dense subset of $C(K)$ are equal; so $\fs\mu=\fs\nu\in\fs(M_p(X))$. The inclusion $\fs(M_p(X))\subset\fs(M_q(X))$ is trivial.

(v)$\Rightarrow$(iv). Let $\mu\in M_q(X)$ and $(h_n)_{n\in\om}\subset C(\be X)$. Put $h'_n=(h_n+c_n)/(2c_n)$, where $c_n=\|h_n\|+1$; then $0\le h'_n\le1$. Let $f\:X\to\Hc$, $f(x)=(h'_n(x))_{n\in\om}$. Then $\fb=(h'_n)_{n\in\om}$, i.e.\ $h'_n=\pi_n\circ\fb$, where $\pi_n\:\Hc\to[0,1]$ is the $n$-th coordinate. By (v) with $K=\Hc$ there is $\nu\in M_p(X)$ with $\fs\nu=\fs\mu$; hence $\nu(h'_n)=\fs\nu(\pi_n)=\fs\mu(\pi_n)=\mu(h'_n)$ for all $n$ and $\nu(1)=\fs\nu(1)=\fs\mu(1)=\mu(1)$, so $\nu(h_n)=\mu(h_n)$ for all $n$.
\end{proof}

The proof of (v)$\Rightarrow$(iv) uses only $K=\Hc$. Hence it suffices to verify (v) for $K=\Hc$; the same remark applies to Proposition~\ref{p:MP} and Theorem~\ref{t:dbR} below.

\begin{note}\label{r:RS}
In \cite[Theorem 2(v)]{RS}, condition (v) was formulated as
\[
M(\be f)(M_q(X))\subset M_p(Y) \tag{$*$}
\]
for every continuous surjection $f\:X\to Y$ onto a separable metrizable space $Y$.
Condition $(*)$ is equivalent to Theorem~\ref{t:main}(v) whenever $M(\be f)$ maps $M_p(X)$ onto $M_p(Y)$ for all such $f$, which is the case, e.g., for $p=a$. In general $(*)$ does not imply (v): if $X$ is pseudocompact, then $Y$ is compact, $M_p(Y)=M(Y)$ for $p\in\{\be,R,\tau,\si,a\}$, and $(*)$ holds automatically for these $p$. Example~\ref{ex:psi} below shows that in this case $(*)$ need not imply (i)--(iv); this is why condition (v) is stated here in terms of images of measures.
\end{note}

\begin{proposition}\label{p:mono}
Let $X$ be a space and $p\pq p'\pq q'\pq q$.
\begin{enumerate}
\item
If $M_p(X)$ is $C$-embedded in $M_q(X)$, then $M_{p'}(X)$ is $C$-embedded in $M_{q'}(X)$.
\item
If $p\pq r\pq q$, $M_p(X)$ is $C$-embedded in $M_r(X)$ and $M_r(X)$ is $C$-embedded in $M_q(X)$, then $M_p(X)$ is $C$-embedded in $M_q(X)$.
\end{enumerate}
\end{proposition}
\begin{proof}
By Theorem~\ref{t:main} we may speak of $G_\de$-density. (i): every nonempty $G_\de$-subset of $M_{q'}(X)$ is the trace of a $G_\de$-subset of $M_q(X)$, which therefore meets $M_p(X)\subset M_{p'}(X)$. (ii): a nonempty $G_\de$-set $H$ of $M_q(X)$ meets $M_r(X)$, and $H\cap M_r(X)$ is a nonempty $G_\de$-subset of $M_r(X)$, so it meets $M_p(X)$.
\end{proof}

\subsection{Probability measures}

\begin{lemma}[positive lifting]\label{l:lift}
Let $p\in\{d,\be,R,\tau,\si,a\}$, let $K$ be a metrizable compact space, $f\:X\to K$ continuous, $\mu\in M^+(\be X)$ and $\nu\in M_p(X)$ with $\fs\nu=\fs\mu$. Then there is $\nu'\in M^+_p(X)$ with $\nu'\le\nu^+$ and $\fs\nu'=\fs\mu$; in particular $\|\nu'\|=\|\mu\|$.
\end{lemma}
\begin{proof}
Put $\la=\fs\mu\ge 0$ and $\rho=\fs\nu^+\ge0$. Since $\fs\nu=\fs\nu^+-\fs\nu^-$, we have $\rho=\la+\fs\nu^-\ge\la$. By the Radon--Nikodym theorem there is a Borel function $w\:K\to[0,1]$ with $\la=w\cdot\rho$. Let $\nu'=(w\circ\fb)\cdot\nu^+$. Then $\nu'\in M^+(\be X)$, $\nu'\le\nu^+\le|\nu|$, so $\nu'\in M_p(X)$ by Lemma~\ref{l:Mp}(i). For $g\in C(K)$,
\[
\fs\nu'(g)=\int (g\circ\fb)(w\circ\fb)\,d\nu^+=\int gw\,d(\fs\nu^+)=\int gw\,d\rho=\la(g),
\]
so $\fs\nu'=\la$. Finally $\|\nu'\|=\nu'(1)=\la(1)=\mu(1)=\|\mu\|$.
\end{proof}

\begin{proposition}\label{p:MP}
Let $X$ be a space and $p\pq q$. The following conditions are equivalent:
\begin{enumerate}
\item
$M_p(X)$ is $C$-embedded in $M_q(X)$;
\item
$P_p(X)$ is $G_\de$-dense in $P_q(X)$;
\item
for every metrizable compact space $K$ and every continuous map $f\:X\to K$,
\[
\fs(P_q(X))\subset\fs(P_p(X)).
\]
\end{enumerate}
\end{proposition}
\begin{proof}
(ii)$\Leftrightarrow$(iii). The proof of the equivalences (iii)$\Leftrightarrow$(iv)$\Leftrightarrow$(v) in Theorem~\ref{t:main} uses only the fact that $M_p(X)\subset M_q(X)\subset M(\be X)$ carry the weak$^*$ topology; it applies verbatim with $M_p(X),M_q(X)$ replaced by $P_p(X),P_q(X)$.

(i)$\Rightarrow$(iii). Let $\mu\in P_q(X)$. By Theorem~\ref{t:main}(v) there is $\nu\in M_p(X)$ with $\fs\nu=\fs\mu$. By Lemma~\ref{l:lift} there is $\nu'\in M^+_p(X)$ with $\fs\nu'=\fs\mu$ and $\|\nu'\|=1$, i.e.\ $\nu'\in P_p(X)$.

(iii)$\Rightarrow$(i). Let $\mu\in M_q(X)$. By Lemma~\ref{l:Mp}(i), $\mu^+,\mu^-\in M^+_q(X)$. If $\mu^+\ne0$, then $\mu^+/\|\mu^+\|\in P_q(X)$, and by (iii) there is $\nu_+\in P_p(X)$ with $\fs\nu_+=\fs\mu^+/\|\mu^+\|$; if $\mu^+=0$, put $\nu_+=0$. Define $\nu_-$ in the same way. Then $\nu=\|\mu^+\|\nu_+-\|\mu^-\|\nu_-\in M_p(X)$ and $\fs\nu=\fs\mu$. By Theorem~\ref{t:main}, $M_p(X)$ is $C$-embedded in $M_q(X)$.
\end{proof}

A dense $C$-embedded subspace of any space is $G_\de$-dense in it \cite[Theorem 6.1.4]{at2009}. Hence, by Proposition~\ref{p:MP}, every positive result on $C$-embeddings of spaces of probability measures yields a result for signed measures. For example, by \cite[Corollaries 5.2 and 5.3]{BanakhChigogidzeFedorchuk2003}, $P_R(Y)$ is $C$-embedded in $P_\si(Y)$ if $Y=X^\kappa$, where $X$ is universally measurable \textup(i.e.\ $M_\si(X)=M_R(X)$\textup), all finite powers of $X$ have countable Souslin number, and $\kappa$ is any cardinal, or if $Y$ is an $AE(0)$-space. Hence $M_R(Y)$ is $C$-embedded in $M_\si(Y)$ for all such $Y$.

\begin{note}\label{r:RS2}
The converse fails: the $C$-embedding of $M_p(X)$ in $M_q(X)$ does not imply the $C^*$-embedding of $P_p(X)$ in $P_q(X)$ (Example~\ref{ex:mrowka}). In the companion paper \cite{RS2} it is shown that $G_\de$-density of $Q_p(X)$ in $Q_q(X)$ is the same condition for all $Q\in\{M,M^+,U,P\}$, and that $C^*$-embedding of $Q_p(X)$ in $Q_q(X)$ implies $G_\de$-density; the spaces $M^+_p(X)$, $U_p(X)$ and $P_p(X)$ are studied there in detail.
\end{note}

%======================================================================
\section{Discrete, compactly supported and Radon measures}\label{sec:dbR}

For a space $X$ let $\cK_d(X)$, $\cK_\be(X)$ and $\cK_R(X)$ be the families of countable, compact and $\si$-compact subsets of $X$, respectively. Thus for $p\in\{d,\be,R\}$ a measure $\mu\in M(\be X)$ belongs to $M_p(X)$ iff it is concentrated on some $S\in\cK_p(X)$. For a continuous map $f$ of $X$ into a metrizable compact space $K$ and $S\in\cK_p(X)$ the set $f(S)$ is countable, compact or $\si$-compact, respectively, hence Borel in $K$.

\begin{lemma}[lifting]\label{l:liftS}
Let $p\in\{d,\be,R\}$, let $K$ be a metrizable compact space, $f\:X\to K$ continuous, $S\in\cK_p(X)$, and let $\la\in M^+(K)$ be concentrated on $f(S)$. Then there is $\nu\in M^+_p(X)$ concentrated on $S$ with $\fs\nu=\la$.
\end{lemma}
\begin{proof}
If $p=d$, then $\la=\sum_{y\in f(S)}\la(\{y\})\de_y$; choose $x_y\in S$ with $f(x_y)=y$ and put $\nu=\sum_y\la(\{y\})\de_{x_y}$. If $p=\be$, apply Lemma~\ref{l:onto}(ii) to the surjection $f|_S\:S\to f(S)$ and the measure $\la|_{f(S)}\in M^+(f(S))$. If $p=R$, write $S=\bigcup_nC_n$ with $C_n$ compact, put $D_n=f(C_n)\setminus\bigcup_{k<n}f(C_k)$, and use the case $p=\be$ to find $\nu_n\in M^+(C_n)$ with $(f|_{C_n})_\sharp\nu_n=\la|_{D_n}$; then $\|\nu_n\|=\la(D_n)$. The Borel sets $D_n$ are pairwise disjoint and $\bigcup_nD_n=f(S)$, so $\la=\sum_n\la|_{D_n}$ and $\sum_n\|\nu_n\|=\la(f(S))<\infty$. Hence $\nu=\sum_n\nu_n$ is a well-defined measure in $M^+(\be X)$, it is concentrated on $S$, and $\fs\nu=\sum_n\la|_{D_n}=\la$.
\end{proof}

\begin{theorem}\label{t:dbR}
Let $X$ be a space, $p\in\{d,\be,R\}$ and $p\pq q$. Then $M_p(X)$ is $C$-embedded \textup($C^*$-embedded\textup) in $M_q(X)$ if and only if for every metrizable compact space $K$, every continuous $f\:X\to K$ and every $\mu\in P_q(X)$ there is $S\in\cK_p(X)$ with
\[
\mu\big(\fb^{-1}(f(S))\big)=1.
\]
\end{theorem}
\begin{proof}
By Theorem~\ref{t:main} and Proposition~\ref{p:MP}, $M_p(X)$ is $C$-embedded (and $C^*$-embedded) in $M_q(X)$ iff $\fs(P_q(X))\subset\fs(P_p(X))$ for every $f$.

Suppose this inclusion holds, and let $\mu\in P_q(X)$. Take $\nu\in P_p(X)$ with $\fs\nu=\fs\mu$, and let $\nu$ be concentrated on $S\in\cK_p(X)$. Then $\mu(\fb^{-1}(f(S)))=\fs\mu(f(S))=\fs\nu(f(S))=\nu(\fb^{-1}(f(S)))\ge\nu(S)=1$.

Conversely, let $\mu\in P_q(X)$ and $\mu(\fb^{-1}(f(S)))=1$ with $S\in\cK_p(X)$. Then $\la=\fs\mu$ is concentrated on $f(S)$, and Lemma~\ref{l:liftS} gives $\nu\in M^+_p(X)$ with $\fs\nu=\la$; since $\|\nu\|=\la(K)=1$, we have $\nu\in P_p(X)$.
\end{proof}

The following sufficient condition is stated in terms of $X$ alone.

\begin{proposition}\label{p:S}
Suppose that for every continuous map $f\:X\to\Hc$ there is a $\si$-compact set $S\subset X$ with $f(S)=f(X)$. Then $M_R(X)$ is $C$-embedded in $M_\si(X)$, and hence $M_\tau(X)$ is $C$-embedded in $M_\si(X)$. If, moreover, $X$ is pseudocompact, then $M_R(X)$ is $C$-embedded in $M(\be X)$.
\end{proposition}
\begin{proof}
Let $f\:X\to\Hc$ be continuous, let $S\subset X$ be a $\si$-compact set with $f(S)=f(X)$, and let $\mu\in P_\si(X)$. If $F\subset\Hc\setminus f(S)$ is compact, then $\fb^{-1}(F)$ is a zero set of $\be X$ disjoint from $X$, so $\fs\mu(F)=\mu(\fb^{-1}(F))=0$ by $(M_\si)$. By inner regularity of $\fs\mu$ on the Borel set $\Hc\setminus f(S)$, we get $\mu(\fb^{-1}(f(S)))=\fs\mu(f(S))=1$. By Theorem~\ref{t:dbR} (with $K=\Hc$, see the remark after Theorem~\ref{t:main}), $M_R(X)$ is $C$-embedded in $M_\si(X)$. The other assertions follow from Proposition~\ref{p:mono} and Corollary~\ref{c:pc}(ii).
\end{proof}

The hypothesis of Proposition~\ref{p:S} holds, for example, if $X$ is pseudo-$\om$-bounded (then $S$ can be taken compact, see the proof of Theorem~\ref{t:beta}) or if every continuous image of $X$ in $\Hc$ is countable (then $S$ can be taken countable); the latter is the case if $\be X$ is scattered (see the proof of Theorem~\ref{t:d}(ii)). It fails for $X=\R^{\mathfrak c}$, although $M_R(X)$ is $C$-embedded in $M_\si(X)$ (Remark~\ref{r:Rc}): for a countable infinite set $J$ of coordinates, the image of $X$ under the projection onto $\R^J\cong(0,1)^\om\subset\Hc$ is not $\si$-compact, since $(0,1)^\om$ is a Baire space whose compact subsets are nowhere dense.

A space $X$ condenses onto a space $Y$ if there is a continuous bijection $X\to Y$.

\begin{cor}\label{c:condense}
Let $X$ be a space that condenses onto a subspace of a separable metrizable space \textup(for example, $X$ is separable metrizable\textup), and let $p\in\{d,\be,R\}$. Then $M_p(X)$ is $C$-embedded \textup($C^*$-embedded\textup) in $M_\tau(X)$ if and only if $M_p(X)=M_\tau(X)$.
\end{cor}
\begin{proof}
Since every separable metrizable space embeds in $\Hc$, there is a continuous $f\:X\to\Hc$ which is injective on $X$. Let $\mu\in P_\tau(X)$, and let $S\in\cK_p(X)$ be as in Theorem~\ref{t:dbR}. The Borel set $B=\fb^{-1}(f(S))\setminus S$ does not meet $X$, because $f$ is injective on $X$. By $(M_\tau)$, $\mu(K)=0$ for every compact $K\subset B$, and by inner regularity of $\mu$ we get $\mu(B)=0$. Hence $\mu(S)=1$ and $\mu\in P_p(X)$. Thus $P_\tau(X)\subset P_p(X)$, and $M_\tau(X)\subset M_p(X)$ by Lemma~\ref{l:Mp}(i).
\end{proof}

\begin{cor}\label{c:classical}
Let $X$ be a separable metrizable space \textup(for example, $X\subset\R$\textup). Then
\begin{enumerate}
\item
$M_R(X)$ is $C$-embedded in $M_\si(X)$ iff $M_R(X)=M_\si(X)$, i.e.\ iff $X$ is universally measurable;
\item
$M_d(X)$ is $C$-embedded in $M_\si(X)$ iff $M_d(X)=M_\si(X)$, i.e.\ iff $X$ has universal measure zero.
\end{enumerate}
\end{cor}
\begin{proof}
A separable metrizable space is Lindel\"of, so $M_\si(X)=M_\tau(X)$ \cite[Corollary 2.2]{BanakhChigogidzeFedorchuk2003}. Apply Corollary~\ref{c:condense}; for the classical meaning of the equalities see \cite{RS}.
\end{proof}

Thus, for subsets of the real line, the functional notions of universal $M_\si$-measurability and universal $M_\si$-measure zero coincide with the classical ones.

\begin{note}\label{r:tau}
(i) By Proposition~\ref{p:mono}, $M_R(X)$ is $C$-embedded in $M_\si(X)$ iff $M_R(X)$ is $C$-embedded in $M_\tau(X)$ and $M_\tau(X)$ is $C$-embedded in $M_\si(X)$.

(ii) For $p=\tau$, Proposition~\ref{p:MP} gives the criterion: $M_\tau(X)$ is $C$-embedded in $M_q(X)$ iff $\fs(P_q(X))\subset\fs(P_\tau(X))$ for every continuous map $f$ of $X$ into a metrizable compact space. We do not know an analogue of Theorem~\ref{t:dbR} for $p=\tau$; the difficulty is that, unlike Radon measures, $\tau$-additive measures do not lift along continuous surjections. For instance, if $X$ is a discrete space of cardinality $\mathfrak c$ and $f\:X\to[0,1]$ is a bijection, then $P_\tau(X)=P_d(X)$ (a measure $\mu\in P(\be X)$ vanishing on the compact set $\be X\setminus X$ is concentrated on a countable subset of the open set $X$), while $P_\tau([0,1])$ contains the Lebesgue measure.
\end{note}

%======================================================================
\section{Embeddings into $M_a(X)=M(\be X)$}\label{sec:a}

\begin{proposition}\label{p:Gdclosed}
For every space $X$ the subspace $M_\si(X)$ is $G_\de$-closed in $M(\be X)$: every $\mu\in M(\be X)\setminus M_\si(X)$ is contained in a $G_\de$-subset of $M(\be X)$ disjoint from $M_\si(X)$.
\end{proposition}
\begin{proof}
Let $\mu\notin M_\si(X)$. There is a zero set $Z=h^{-1}(0)$, $h\in C(\be X)$, $0\le h\le1$, with $Z\cap X=\es$ and $|\mu|(Z)>0$. The restriction $\mu|_Z$ is a nonzero Radon measure, so there is $g\in C(\be X)$ with $\int_Z g\,d\mu\ne0$. Let $\psi_n=\max\{0,1-nh\}\in C(\be X)$; then $0\le\psi_n\le1$ and $\psi_n\downarrow\mathbf 1_Z$ pointwise. Put $H=\set{\nu\in M(\be X):\nu(g\psi_n)=\mu(g\psi_n)\text{ for all }n}$, a $G_\de$-set containing $\mu$. For $\nu\in H$, by dominated convergence, $\int_Zg\,d\nu=\lim_n\nu(g\psi_n)=\int_Zg\,d\mu\ne0$, so $|\nu|(Z)>0$ and $\nu\notin M_\si(X)$.
\end{proof}

\begin{cor}\label{c:pc}
Let $X$ be a space.
\begin{enumerate}
\item
$M_\si(X)=M(\be X)$ iff $M_\si(X)$ is $C$-embedded \textup($C^*$-embedded\textup) in $M(\be X)$ iff $X$ is pseudocompact.
\item
Let $p\pq\si$. Then $M_p(X)$ is $C$-embedded in $M(\be X)$ iff $X$ is pseudocompact and $M_p(X)$ is $C$-embedded in $M_\si(X)$.
\end{enumerate}
\end{cor}
\begin{proof}
(i) If $X$ is pseudocompact, the condition $(M_\si)$ is void by Lemma~\ref{l:pc}. If $X$ is not pseudocompact, there is a nonempty zero set $Z\subset\be X\setminus X$ (Lemma~\ref{l:pc}), and $\de_z\notin M_\si(X)$ for $z\in Z$. By Proposition~\ref{p:Gdclosed}, $M_\si(X)$ is $G_\de$-dense in $M(\be X)$ only if $M_\si(X)=M(\be X)$. It remains to apply Theorem~\ref{t:main}.

(ii) If $M_p(X)$ is $C$-embedded in $M(\be X)$, then so is $M_\si(X)$ (Proposition~\ref{p:mono}); by (i), $X$ is pseudocompact and $M_\si(X)=M(\be X)$. The converse is obvious from (i).
\end{proof}

\begin{theorem}\label{t:beta}
For a space $X$ the following are equivalent:
\begin{enumerate}
\item
$X$ is pseudo-$\om$-bounded;
\item
$M_\be(X)$ is $C$-embedded in $M(\be X)$;
\item
$M_\be(X)$ is $C^*$-embedded in $M_q(X)$ for some $q\in\{R,\tau,\si,a\}$.
\end{enumerate}
Consequently, if $X$ is pseudo-$\om$-bounded, then $M_p(X)$ is $C$-embedded in $M_q(X)$ for all $\be\pq p\pq q$.
\end{theorem}
\begin{proof}
(i)$\Rightarrow$(ii). We verify Theorem~\ref{t:main}(v) for $p=\be$, $q=a$. Let $f\:X\to\Hc$ be continuous and $K=\fb(\be X)=\cl{f(X)}$. Let $\{B_n\}$ be a countable base of nonempty open subsets of $K$. Since $f(X)$ is dense in $K$, the sets $U_n=f^{-1}(B_n)$ are nonempty and open in $X$. Let $C\subset X$ be a compact set meeting all $U_n$. Then $f(C)$ is a compact subset of $K$ meeting every $B_n$, hence $f(C)=K$. By Lemma~\ref{l:onto}(ii) applied to $f|_C\:C\to K$, every measure in $M(K)$ is the image of a measure on $C$; measures on $C$ belong to $M_\be(X)$. Since $\fs(M(\be X))\subset M(K)$, we get $\fs(M(\be X))\subset\fs(M_\be(X))$.

(ii)$\Rightarrow$(iii) is trivial. If (iii) holds, then $M_\be(X)$ is $G_\de$-dense in $M_q(X)$ (Theorem~\ref{t:main}), hence in $M_R(X)\subset M_q(X)$, and thus Theorem~\ref{t:main}(v) holds for $p=\be$, $q=R$.

(iii)$\Rightarrow$(i). Let $(U_n)_{n\ge1}$ be nonempty open sets in $X$. Choose $x_n\in U_n$ and continuous $h_n\:X\to[0,1]$ with $h_n(x_n)=1$ and $h_n=0$ outside $U_n$. Let $f=(h_1,h_2,\dots)\:X\to\Hc$ and $y_n=f(x_n)$. The measure $\mu=\sum_n2^{-n}\de_{x_n}$ belongs to $M_d(X)\subset M_R(X)$, and $\fs\mu(\{y_n\})\ge2^{-n}>0$. By Theorem~\ref{t:main}(v) there is $\nu\in M_\be(X)$ with $\fs\nu=\fs\mu$; let $\nu$ be concentrated on a compact $C\subset X$. Then $\nu(\fb^{-1}(y_n))=\fs\nu(\{y_n\})\ne0$, so $\fb^{-1}(y_n)\cap C\ne\es$: there is $c\in C$ with $f(c)=y_n$, whence $h_n(c)=1$ and $c\in U_n$. Thus $C$ meets all $U_n$.

The last assertion follows from (ii) and Proposition~\ref{p:mono}.
\end{proof}

\begin{note}\label{r:izv}
Corollary~\ref{c:pc}(i) and Theorem~\ref{t:beta} are signed counterparts of \cite[Theorems 1.22 and 1.27]{Reznichenko2025betap}: $X$ is pseudocompact iff $P_\si(X)$ is $C^*$-embedded in $P(\be X)$, and $X$ is pseudo-$\om$-bounded iff $P_\be(X)$ is $C^*$-embedded in $P(\be X)$. Thus for $p\in\{\be,\si\}$ the $C$-embedding of $M_p(X)$ in $M(\be X)$ is equivalent to the $C^*$-embedding of $P_p(X)$ in $P(\be X)$. For $p\in\{d,R,\tau\}$ this is not the case (Example~\ref{ex:mrowka}).
\end{note}

\begin{theorem}\label{t:d}
Let $X$ be a space.
\begin{enumerate}
\item
$M_d(X)$ is $C$-embedded \textup($C^*$-embedded\textup) in $M_R(X)$ iff $M_d(X)=M_R(X)$ iff every compact subset of $X$ is scattered.
\item
$M_d(X)$ is $C$-embedded \textup($C^*$-embedded\textup) in $M(\be X)$ iff $\be X$ is scattered.
\end{enumerate}
\end{theorem}
\begin{proof}
(i) If all compact subsets of $X$ are scattered, then every $\mu\in M_R(X)$ is concentrated on a union $\bigcup_nC_n$ of compact scattered sets, the restrictions $\mu|_{C_n}$ are discrete by (RK), and so $\mu\in M_d(X)$. Assume now that $M_d(X)$ is $C^*$-embedded in $M_R(X)$ and $C\subset X$ is a compact set which is not scattered. By (PS) there is a continuous $g\:C\to[0,1]$ onto; extend it to $\bar g\:\be X\to[0,1]$ ($\be X$ is normal) and put $f=\bar g|_X\:X\to[0,1]$. By Lemma~\ref{l:onto}(ii) there is $\mu\in P(C)\subset M_R(X)$ such that $\fs\mu$ is the Lebesgue measure on $[0,1]$. By Theorem~\ref{t:main}(v), $\fs\mu=\fs\nu$ for some $\nu\in M_d(X)$; but $\fs\nu$ is discrete, a contradiction.

(ii) Suppose $M_d(X)$ is $C^*$-embedded in $M(\be X)$. Then $X$ is pseudocompact by Corollary~\ref{c:pc}. If $\be X$ is not scattered, take $g\:\be X\to[0,1]$ onto (PS) and repeat the argument of (i) with $C$ replaced by $\be X$ and $f=g|_X$.

Conversely, let $\be X$ be scattered. Then $X$ is pseudocompact by Lemma~\ref{l:scat}. Let $f\:X\to\Hc$ be continuous. Then $K=\fb(\be X)=f(X)$ is a metrizable compact space, and it is countable: otherwise $K$ contains a Cantor set, so by (PS) $K$, and hence $\be X$, admits a continuous map onto $[0,1]$, contradicting (PS) for $\be X$. So every $\la\in M(K)$ has the form $\la=\sum_ia_i\de_{y_i}$ with $y_i\in f(X)$, $\sum|a_i|<\infty$. Choose $x_i\in f^{-1}(y_i)$; then $\nu=\sum_ia_i\de_{x_i}\in M_d(X)$ and $\fs\nu=\la$. Hence $\fs(M(\be X))\subset\fs(M_d(X))$, and Theorem~\ref{t:main} applies.
\end{proof}

\begin{theorem}\label{t:R}
For a space $X$ the following are equivalent:
\begin{enumerate}
\item
$M_R(X)$ is $C$-embedded \textup($C^*$-embedded\textup) in $M(\be X)$;
\item
$P_R(X)$ is $G_\de$-dense in $P(\be X)$;
\item
$X$ is pseudocompact, and for every continuous map $f$ of $X$ onto a metrizable compact space $K$ and every $\la\in P(K)$ there is a Radon probability measure $\nu$ on $X$ with $f_\sharp\nu=\la$.
\end{enumerate}
These conditions hold if $P_R(X)$ is pseudocompact or $X$ is pseudo-$\om$-bounded. If every compact subset of $X$ is scattered, then \textup{(i)--(iii)} are equivalent to the scatteredness of $\be X$.
\end{theorem}
\begin{proof}
(i)$\Leftrightarrow$(ii) is Proposition~\ref{p:MP}.

(i)$\Rightarrow$(iii). $X$ is pseudocompact by Corollary~\ref{c:pc}. Let $f\:X\to K$ be onto. Since $X$ is pseudocompact, $\fb(\be X)=K$, and by Lemma~\ref{l:onto}(ii) every $\la\in P(K)$ is $\fs\mu$ for some $\mu\in P(\be X)$. By Proposition~\ref{p:MP}, $\la=\fs\nu$ for some $\nu\in P_R(X)$.

(iii)$\Rightarrow$(ii). For a continuous map $f$ of $X$ into a metrizable compact space the set $K=f(X)$ is compact, and $\fs(P(\be X))\subset P(K)\subset\fs(P_R(X))$ by (iii). Apply Proposition~\ref{p:MP}.

If $P_R(X)$ is pseudocompact and $f\:X\to K$ is a continuous map into a metrizable compact space, then $\fs(P_R(X))$ is a pseudocompact subset of the metrizable space $P(K)$, hence compact. It contains all finitely supported probability measures on $f(X)$, which are dense in $P(\cl{f(X)})$ by Lemma~\ref{l:onto}(i). Hence $\fs(P(\be X))\subset P(\cl{f(X)})=\fs(P_R(X))$, and (ii) holds by Proposition~\ref{p:MP}. If $X$ is pseudo-$\om$-bounded, (i) follows from Theorem~\ref{t:beta}. If all compact subsets of $X$ are scattered, then $M_R(X)=M_d(X)$ by Theorem~\ref{t:d}(i), and Theorem~\ref{t:d}(ii) applies.
\end{proof}

\begin{note}\label{r:R}
In \cite[Theorem 1.9]{Reznichenko2025betap} it is shown that $P_R(X)$ is $C^*$-embedded in $P(\be X)$ iff $P_R(X)$ is pseudocompact; by Theorem~\ref{t:R}, this condition implies (i)--(iii). The converse is false (Example~\ref{ex:mrowka}).

By Theorem~\ref{t:dbR} and Lemma~\ref{l:onto}(ii), condition (iii) of Theorem~\ref{t:R} can be replaced by: $X$ is pseudocompact, and for every continuous map $f$ of $X$ onto a metrizable compact space $K$, every $\la\in P(K)$ is concentrated on $f(S)$ for some $\si$-compact $S\subset X$. Equivalently: $X$ is pseudocompact, and for every continuous map $f\:X\to\Hc$ no nonzero measure $\la\in M^+(f(X))$ vanishes on all sets $f(C)$ with $C\subset X$ compact. Indeed, if some $\la\in P(f(X))$ is not concentrated on any set $f(S)$, choose a $\si$-compact $S_0\subset X$ with $\la(f(S_0))$ maximal (a countable union of $\si$-compact sets is $\si$-compact); then $\la|_{f(X)\setminus f(S_0)}$ is a nonzero measure vanishing on all sets $f(C)$. The converse implication is obvious.

By Proposition~\ref{p:S}, conditions (i)--(iii) hold if $X$ is pseudocompact and every continuous image $f(X)\subset\Hc$ of $X$ is the image $f(S)$ of a $\si$-compact set $S\subset X$; this covers pseudo-$\om$-bounded spaces and spaces with scattered $\be X$. This sufficient condition is not necessary (Example~\ref{ex:83}). The set $S$ in Theorem~\ref{t:dbR} cannot, in general, be chosen independently of $f$ (Remark~\ref{r:Sdep}), nor, for a fixed $f$, independently of $\la$ (Remark~\ref{r:83}).
\end{note}

%======================================================================
\section{Examples}\label{sec:example}

Recall the construction of Mr\'owka--Isbell spaces \cite{Mrowka1954}, \cite[3.6.I]{EngelkingBookRu}. A family $\cA$ of infinite subsets of $\N$ is almost disjoint if $A\cap B$ is finite for distinct $A,B\in\cA$, and maximal almost disjoint (MAD) if, moreover, every infinite $B\subset\N$ meets some $A\in\cA$ in an infinite set. The space $\Psi(\cA)=\N\cup\cA$ has the topology in which the points of $\N$ are isolated and a base at $A\in\cA$ is formed by the sets $\{A\}\cup(A\setminus F)$, $F$ finite. $\Psi(\cA)$ is a locally compact Tychonoff space, $\N$ is dense in it, and $\cA$ is closed and discrete.

\begin{example}\label{ex:psi}
There is a pseudocompact locally compact space $X$ such that
\begin{enumerate}
\item
$M_d(X)=M_R(X)=M_\tau(X)$ and $M_\si(X)=M(\be X)$;
\item
$M_\tau(X)$ is not $C^*$-embedded in $M_\si(X)$, and $M_R(X)$ is not $C^*$-embedded in $M_a(X)=M(\be X)$.
\end{enumerate}
\end{example}
\begin{proof}
Let $(r_n)_{n\in\N}$ be an enumeration of the rationals in $[0,1]$. Call an infinite set $A\subset\N$ convergent if the sequence $(r_n)_{n\in A}$ converges in $[0,1]$. By Zorn's lemma there is an almost disjoint family $\cA$ of convergent sets which is maximal among such families. Then $\cA$ is MAD: if $B\subset\N$ is infinite, it contains an infinite convergent subset $B'$ (by compactness of $[0,1]$), and $B'$ meets some $A\in\cA$ in an infinite set by maximality; so does $B$.

Let $X=\Psi(\cA)$.

First, $X$ is pseudocompact. Indeed, let $g\:X\to\R$ be continuous and unbounded. Since $\N$ is dense, $g$ is unbounded on $\N$, so there is an infinite $B\subset\N$ with $|g(n)|\to\infty$ along $B$. Take $A\in\cA$ with $A\cap B$ infinite. Every neighbourhood of $A$ contains all but finitely many points of $A\cap B$, which contradicts continuity of $g$ at $A$.

Second, compact subsets of $X$ are countable: a compact $C\subset X$ meets the closed discrete set $\cA$ in a finite set, so $C\subset\N\cup(C\cap\cA)$.

Third, $X$ admits a continuous map onto $[0,1]$. Indeed, let $h(n)=r_n$ for $n\in\N$ and $h(A)=\lim_{n\in A}r_n$ for $A\in\cA$. Then $h$ is continuous at each $A$ by the definition of the topology, and $h(X)$ is compact (as $X$ is pseudocompact) and contains $\mathbb Q\cap[0,1]$, so $h(X)=[0,1]$. In particular $\be X$ is not scattered.

(i) Since $X$ is locally compact, $X$ is open in $\be X$ and $\be X\setminus X$ is compact. Hence $\mu\in M_\tau(X)$ implies $|\mu|(\be X\setminus X)=0$; by inner regularity of $|\mu|$, $\mu$ is concentrated on a $\si$-compact subset of $X$, which is countable. So $M_\tau(X)\subset M_d(X)$, and the converse inclusions are obvious. $M_\si(X)=M(\be X)$ by Corollary~\ref{c:pc}(i).

(ii) By Theorem~\ref{t:R} (the case of scattered compact subsets), $M_R(X)$ is not $C^*$-embedded in $M(\be X)$. Since $M_\tau(X)=M_R(X)$ and $M_\si(X)=M(\be X)$, $M_\tau(X)$ is not $C^*$-embedded in $M_\si(X)$.

For the reader's convenience, here is a direct argument for (ii). Let $\bar h\:\be X\to[0,1]$ extend $h$, and let $\mu\in P(\be X)$ satisfy $\bar h_\sharp\mu=\la$, the Lebesgue measure (Lemma~\ref{l:onto}(ii)). Let $(g_n)$ be dense in $C[0,1]$. The $G_\de$-set $H=\set{\nu\in M(\be X):\nu(g_n\circ\bar h)=\mu(g_n\circ\bar h)\ \forall n}$ contains $\mu$, and every $\nu\in H$ satisfies $\bar h_\sharp\nu=\la$. If $\nu\in M_\tau(X)=M_d(X)$, then $\bar h_\sharp\nu$ is discrete, which is impossible. Hence $H\cap M_\tau(X)=\es$, and Theorem~\ref{t:main} applies.
\end{proof}

\begin{note}
Example~\ref{ex:psi} shows that Corollaries~1 and~2 of \cite{RS} require additional assumptions on $X$: not every pseudocompact space is universally $M^*_{[R,a]}$-measurable (universally $M_a^*$-measurable in the terminology of \cite{RS}), and not every space is universally $M_{[\tau,\si]}$-measurable. By Proposition~\ref{p:MP}, for this $X$ the space $P_\tau(X)$ is not $G_\de$-dense, and hence not $C$-embedded, in $P_\si(X)=P(\be X)$. So for this space the answers to the questions (BCF) and (B1) are negative both for probability and for signed measures. The space $X$ above has a non-scattered Stone--\v Cech compactification; this is exactly what distinguishes it from the Mr\'owka space of \cite[Example 1.15]{Reznichenko2025betap}, for which the signed versions do hold, as the next example shows.
\end{note}

\begin{example}\label{ex:mrowka}
Let $X$ be the Mr\'owka space of \cite[Example 1.15, Theorem 5.6]{Reznichenko2025betap}: $X=\Psi(D,\mathcal E)$ is a Mr\'owka--Isbell space with $|\be X\setminus X|=1$ such that $P_R(X)$ is not pseudocompact \textup(see also \cite{Mrowka1977}, \cite[Section 8.6]{Hernandez-HernandezHrusak2018}\textup). Here $D$ is a countable set of isolated points, $\mathcal E=\set{E_P:P\in\mathcal P}$ is an almost disjoint family of infinite subsets of $D$ indexed by a set $\mathcal P$ disjoint from $D$, the points $P\in\mathcal P$ are the nonisolated points of $X$, and $\{P\}\cup E_P$ is a compact open neighbourhood of $P$. Then $M_d(X)=M_R(X)=M_\tau(X)$ is $C$-embedded in $M_\si(X)=M(\be X)$, while $P(X)=P_\tau(X)$ is not pseudocompact.
\end{example}
\begin{proof}
As in Example~\ref{ex:psi}, $X$ is pseudocompact and locally compact, its compact subsets are countable, $M_d(X)=M_R(X)=M_\tau(X)$, and $M_\si(X)=M(\be X)$. The space $\be X=X\cup\{\infty\}$ is scattered: a nonempty $S\subset\be X$ has an isolated point, namely any point of $S\cap D$ if $S\cap D\ne\es$; otherwise any point $P\in S\setminus\{\infty\}$, which is isolated in $\be X\setminus D$ because $\{P\}\cup E_P$ is a neighbourhood of $P$ in $\be X$; and $\infty$ if $S=\{\infty\}$. By Theorem~\ref{t:R} (the case of scattered compact subsets), $M_R(X)$ is $C$-embedded in $M(\be X)$. The space $P(X)$ is not pseudocompact by \cite[Example 1.15]{Reznichenko2025betap}.
\end{proof}

\begin{note}
By \cite[Theorem 1.9]{Reznichenko2025betap}, for the space $X$ of Example~\ref{ex:mrowka} the space $P(X)=P_\tau(X)$ is not $C^*$-embedded in $P(\be X)=P_\si(X)$, although it is $G_\de$-dense in it (Proposition~\ref{p:MP}). This confirms the observation of \cite{RS} that the passage from probability to signed measures can change the answers to the questions (BCF) and (B1); Example~\ref{ex:psi} shows that this depends on the space. It also shows that the $C$-embedding of $M_R(X)$ in $M(\be X)$ does not imply the pseudocompactness of $P_R(X)$. The spaces $M^+_R(X)$ and $U_R(X)$ for this $X$ are considered in \cite{RS2}.
\end{note}

\begin{example}\label{ex:sigma}
Let $K=[0,1]^{\om_1}$ and let $\Sigma=\set{x\in K:|\{\alpha<\om_1:x_\alpha\ne0\}|\le\om}$ be the $\Sigma$-product. Then
\begin{enumerate}
\item
$\Sigma$ is $\om$-bounded and $\be\Sigma=K$;
\item
$M_R(\Sigma)\ne M_\tau(\Sigma)\ne M_\si(\Sigma)=M(K)$;
\item
$M_\be(\Sigma)$ is $C$-embedded in $M(K)$; consequently, each of the spaces $M_R(\Sigma)\subset M_\tau(\Sigma)\subset M_\si(\Sigma)$ is $C$-embedded in each larger one.
\end{enumerate}
Thus $C$-embeddings $M_R(X)\subset M_\tau(X)\subset M_\si(X)$ may hold with all inclusions proper.
\end{example}
\begin{proof}
For $J\subset\om_1$ let $\pi_J\:K\to[0,1]^J$ be the projection. We use two standard facts: every $G_\de$-subset $H$ of $K$ and every Baire subset $B$ of $K$ is determined by countably many coordinates, i.e.\ for every $y\in H$ (resp.\ for $B$) there is a countable $J$ with $\pi_J^{-1}(\pi_J(y))\subset H$ (resp.\ $B=\pi_J^{-1}(\pi_J(B))$).

(i) If $A\subset\Sigma$ is countable, the union $J$ of the supports of its points is countable, and $\cl A\subset[0,1]^J\times\{0\}^{\om_1\setminus J}\subset\Sigma$, a compact set. So $\Sigma$ is $\om$-bounded, hence pseudocompact. If $H$ is a nonempty $G_\de$-set in $K$, $y\in H$ and $\pi_J^{-1}(\pi_J(y))\subset H$ with $J$ countable, then the point equal to $y$ on $J$ and to $0$ outside $J$ lies in $H\cap\Sigma$. So $\Sigma$ is $G_\de$-dense in $K$, hence $C$-embedded in $K$ by Theorem~\ref{t:prod:main}, and $\be\Sigma=K$.

(ii) Since $\Sigma$ is pseudocompact, $M_\si(\Sigma)=M(K)$ (Corollary~\ref{c:pc}). For $p\in K\setminus\Sigma$ the compact set $\{p\}$ misses $\Sigma$, so $\de_p\notin M_\tau(\Sigma)$. Let $\la$ be the product of the Lebesgue measures on $K$. It is completion regular \cite{Kakutani1943}: for every compact $F\subset K$ there are Baire sets $B_1\subset F\subset B_2$ with $\la(B_2\setminus B_1)=0$. If $F\cap\Sigma=\es$ and $\la(F)>0$, then $B_1\neq\es$ and $B_1=\pi_J^{-1}(\pi_J(B_1))$ for a countable $J$, so $B_1$ meets $\Sigma$, a contradiction. Hence $\la(F)=0$ for every compact $F\subset K\setminus\Sigma$, i.e.\ $\la\in P_\tau(\Sigma)$. If $F\subset\Sigma$ is compact and $\la(F)>0$, then $\es\ne B_1=\pi_J^{-1}(\pi_J(B_1))\subset\Sigma$, which is impossible, since every fibre of $\pi_J$ contains points with uncountable support. Hence $\la(S)=0$ for every $\si$-compact $S\subset\Sigma$, and $\la\notin M_R(\Sigma)$.

(iii) follows from (i), Theorem~\ref{t:beta} and Proposition~\ref{p:mono}.
\end{proof}

\begin{note}\label{r:Sdep}
In Theorem~\ref{t:dbR} the set $S$ depends, in general, on $f$, even for $p=R$ and $q=a$. Let $\la$ be the measure from the proof of Example~\ref{ex:sigma}(ii), and let $S\subset\Sigma$ be $\si$-compact. Then $\la(S)=0$; by completion regularity there is a Baire set $B\supset S$ with $\la(B)=0$, and $B=\pi_J^{-1}(\pi_J(B))$ for a countable $J\subset\om_1$. For the map $f=\pi_J|_\Sigma\:\Sigma\to[0,1]^J$ we have $\fb=\pi_J$ (as $\be\Sigma=K$), and hence $\la(\fb^{-1}(f(S)))\le\la(B)=0$. Thus no single $\si$-compact set $S\subset\Sigma$ satisfies the condition of Theorem~\ref{t:dbR} for $\mu=\la$ and all $f$, although $M_R(\Sigma)$ is $C$-embedded in $M(K)$. For a fixed map $f$, the set $S$ must in general depend on $\mu$ (Remark~\ref{r:83}).
\end{note}

\begin{example}\label{ex:tree}
Let $D=2^{<\om}$ be the set of finite binary sequences, and for $x\in2^\om$ let $C_x=\set{x|n:n\in\om}$. The family $\set{C_x:x\in2^\om}$ is almost disjoint. Let $X=D\cup2^\om$ be the corresponding Mr\'owka--Isbell space: the points of $D$ are isolated, and the sets $\{x\}\cup\set{x|n:n\ge k}$, $k\in\om$, form a base at $x\in2^\om$. Then
\begin{enumerate}
\item
$X$ is locally compact, submetrizable and realcompact, and $M_d(X)=M_R(X)=M_\tau(X)$;
\item
$M_\tau(X)$ is not $C^*$-embedded in $M_\si(X)$.
\end{enumerate}
\end{example}
\begin{proof}
Let $K$ be the set $D\cup2^\om$ with the topology of the Cantor tree: $2^\om$ carries its usual topology, the points of $D$ are isolated, and a sequence of nodes converges to $x\in2^\om$ iff the nodes eventually extend $x|k$ for every $k$. Then $K$ is a metrizable compact space and the identity map $f\:X\to K$ is a continuous bijection, so $X$ is submetrizable.

(i) $X$ is locally compact, and its compact subsets are countable, since $2^\om$ is closed and discrete in $X$. As in Example~\ref{ex:psi}, $M_d(X)=M_R(X)=M_\tau(X)$. The space $X$ is realcompact, since it admits a continuous bijection onto the metrizable compact space $K$ \cite[Corollaries 8.15 and 8.18]{GillmanJerison1960}.

(ii) For $h\in C(\be X)$ the restriction $h|_{2^\om}$ is a Borel (even Baire class one) function on $2^\om$, because $h(x)=\lim_nh(x|n)$ and $x\mapsto h(x|n)$ is locally constant. Let $m$ be the Cantor measure on $2^\om$ and define $\mu\in P(\be X)$ by $\mu(h)=\int_{2^\om}h\,dm$. Then $\mu\in P_\si(X)$: if $Z=h^{-1}(0)$, $0\le h\le1$, is a zero set of $\be X$ disjoint from $X$, then $\psi_n=\max\{0,1-nh\}\downarrow\mathbf 1_Z$ and $\psi_n\to0$ on $X$, so $\mu(Z)=\lim_n\mu(\psi_n)=\lim_n\int_{2^\om}\psi_n\,dm=0$. For $g\in C(K)$ we have $\fs\mu(g)=\int_{2^\om}g\,dm$, i.e.\ $\fs\mu=m$ is a measure on $K$ vanishing on points. On the other hand, $\fs\nu$ is discrete for every $\nu\in P_\tau(X)=P_d(X)$. By Proposition~\ref{p:MP}, $M_\tau(X)$ is not $C$-embedded in $M_\si(X)$, and by Theorem~\ref{t:main} it is not $C^*$-embedded in $M_\si(X)$.
\end{proof}

Thus realcompactness, even together with local compactness, does not imply that $M_\tau(X)$ is $C$-embedded in $M_\si(X)$. This is consistent with the fact that a locally compact realcompact space need not be measure-compact, i.e.\ need not satisfy $M_\tau(X)=M_\si(X)$.

\begin{note}\label{r:Rc}
Conversely, measure-compactness is not necessary for the $C$-embedding of $M_\tau(X)$ in $M_\si(X)$, even for realcompact $X$. Let $X=\R^{\mathfrak c}$. Then $X$ is realcompact and $M_\tau(X)\ne M_\si(X)$ \cite{KempermanMaharam1970} \textup(see also \cite[Section 2]{BanakhChigogidzeFedorchuk2003}\textup). On the other hand, $P_R(X)$ is $C$-embedded in $P_\si(X)$ by \cite[Corollary 5.2]{BanakhChigogidzeFedorchuk2003}, so $M_R(X)$, and hence $M_\tau(X)$, is $C$-embedded in $M_\si(X)$ (see the paragraph after Proposition~\ref{p:MP}, and Proposition~\ref{p:mono}).
\end{note}

The next proposition shows that $M_\tau(X)$ is $C$-embedded in $M(\be X)$ for every $G_\de$-dense subspace $X$ of a product of metrizable compact spaces, for instance, for the $\Sigma$-product of Example~\ref{ex:sigma}. Unlike in Example~\ref{ex:sigma}, here $X$ need not contain large compact sets (Example~\ref{ex:cube}).

\begin{proposition}\label{p:tau}
Let $T=\prod_{s\in S}T_s$ be a product of nonempty metrizable compact spaces, and let $X$ be a $G_\de$-dense subspace of $T$. Then $X$ is pseudocompact, $\be X=T$, and $M_\tau(X)$ is $C$-embedded in $M_\si(X)=M(\be X)$.
\end{proposition}
\begin{proof}
For $J\subset S$ let $T_J=\prod_{s\in J}T_s$ and let $\pi_J\:T\to T_J$ be the projection. By Theorem~\ref{t:prod:main}, $X$ is $C$-embedded in $T$, so $\be X=T$, and $X$ is pseudocompact by Lemma~\ref{l:pc}; hence $M_\si(X)=M(T)$ (Corollary~\ref{c:pc}). By Proposition~\ref{p:MP} it suffices to show that $\fs(P(T))\subset\fs(P_\tau(X))$ for every continuous map $f$ of $X$ into a metrizable compact space $L$. The extension $\fb\:T\to L$ has the form $\fb=g\circ\pi_J$ for a countable $J\subset S$ and a continuous $g\:T_J\to L$ \cite[Problem 2.7.12(d)]{EngelkingBookRu}. If $J=S$, then $T$ is metrizable, its points are $G_\de$-sets, and $X=T$; in this case there is nothing to prove. Let $S'=S\setminus J\ne\es$. For $s\in S'$ choose a continuous surjection $\phi_s\:2^\om\to T_s$ \cite[Problem 4.5.9(b)]{EngelkingBookRu}, let $\phi=\prod_{s\in S'}\phi_s\:D\to T_{S'}$, where $D=(2^\om)^{S'}$, let $h$ be the Haar measure of the Cantor cube $D$ \textup(the product of the measures $\frac12(\de_0+\de_1)$\textup), and put $m=\phi_\sharp h\in P(T_{S'})$.

Let $\la=\fs\mu$ with $\mu\in P(T)$. Put $\nu_J=(\pi_J)_\sharp\mu\in P(T_J)$ and $\nu=\nu_J\otimes m\in P(T)$. Then $\fs\nu=g_\sharp(\pi_J)_\sharp\nu=g_\sharp\nu_J=\fb_\sharp\mu=\la$.

We show that $\nu\in P_\tau(X)$. Let $F\subset T\setminus X$ be compact and suppose $\nu(F)>0$. By Fubini's theorem, $m(F_u)>0$ for some $u\in T_J$, where $F_u=\set{v\in T_{S'}:(u,v)\in F}$. So $h(\phi^{-1}(F_u))>0$. The measure $h$ is completion regular \cite{Kakutani1943}, so the compact set $\phi^{-1}(F_u)$ contains a nonempty Baire subset of $D$, and hence a set of the form $\set{e\in D:e|_E=w}$, where $E\subset\om\times S'$ is countable and $w\in2^E$. Let $S''$ be the projection of $E$ to $S'$. The image of the latter set under $\phi$ is $\prod_{s\in S''}Q_s\times T_{S'\setminus S''}$, where the sets $Q_s\subset T_s$ are nonempty. Taking $z\in\prod_{s\in S''}Q_s$, we get $\pi_{J\cup S''}^{-1}(u,z)\subset F$. This is a nonempty $G_\de$-subset of $T$ disjoint from $X$, which contradicts the $G_\de$-density of $X$. Hence $\nu(F)=0$ for every compact $F\subset T\setminus X$, i.e.\ $\nu\in P_\tau(X)$.
\end{proof}

\begin{example}\label{ex:cube}
Let $\kappa=\mathfrak c$, let $K=I^\kappa$, where $I=[0,1]$ or $I=\{0,1\}$, and let $X$ be a pseudocompact dense subspace of $K$ all of whose compact subsets are countable \textup(such spaces exist, see below\textup). Then
\begin{enumerate}
\item
$\be X=K$, $M_R(X)=M_d(X)$, $M_\si(X)=M(K)$ and $M_\tau(X)\ne M_\si(X)$;
\item
$M_\tau(X)$ is $C$-embedded in $M_\si(X)=M(\be X)$;
\item
$M_R(X)$ is not $C^*$-embedded in $M(\be X)$, and $M_R(X)$ is not $C^*$-embedded in $M_\tau(X)$.
\end{enumerate}
\end{example}

Such spaces $X$ exist in ZFC. Shakhmatov \cite{Shakhmatov1986} and the first author \cite{Reznichenko1989} constructed pseudocompact dense subspaces of $[0,1]^{\mathfrak c}$ all of whose countable subsets are closed and discrete. Recently it was shown that pseudocompact Korovin orbits also have this property \cite{ReznichenkoTkachenko2024}; they give further examples, also for $I=\{0,1\}$. Let $G$ be an Abelian group of cardinality $\mathfrak c$, let $f\:G\to I$ be a Korovin mapping \cite{Korovin1992} \textup(see \cite[Theorem 2.3]{ReznichenkoTkachenkoKorovin}\textup), and let $X=Gf\subset I^G$ be the corresponding Korovin orbit. Then $\pi_J(X)=I^J$ for every countable $J\subset G$; hence $X$ is $G_\de$-dense in $I^G$ and therefore pseudocompact \textup(Theorem~\ref{t:prod:main} and Lemma~\ref{l:pc}\textup). Every countable subset of $X$ is closed and discrete \cite[Theorem 2]{ReznichenkoTkachenko2024}, so all compact subsets of $X$ are finite.

\begin{proof}
For $J\subset\kappa$ let $\pi_J\:K\to I^J$ be the projection. (i) For every countable $J\subset\kappa$ the set $\pi_J(X)$ is a pseudocompact dense subset of the metrizable compact space $I^J$, hence compact and equal to $I^J$. By Theorem~\ref{t:prod:main} (applied to $\Pi=Z=K$ and $A=X$), $X$ is $G_\de$-dense in $K$. By Proposition~\ref{p:tau}, $\be X=K$ and $M_\si(X)=M(K)$. Since compact subsets of $X$ are countable, $M_R(X)=M_d(X)$. Finally, $X\ne K$ (as $K$ is an uncountable compact space), and $\de_p\notin M_\tau(X)$ for $p\in K\setminus X$.

(ii) follows from (i) and Proposition~\ref{p:tau}.

(iii) The compact space $K=\be X$ is not scattered (it maps onto $[0,1]$), so by Theorem~\ref{t:d}(ii) the space $M_d(X)=M_R(X)$ is not $C^*$-embedded in $M(\be X)$. If $M_R(X)$ were $C^*$-embedded in $M_\tau(X)$, then by (ii) and Proposition~\ref{p:mono}(ii) it would be $C$-embedded in $M(\be X)$.
\end{proof}

Thus for a space $X$ as in Example~\ref{ex:cube}, $M_\tau(X)$ is $C$-embedded in $M(\be X)$ while $M_R(X)$ is not, and $M_\tau(X)\ne M_\si(X)$, $M_\tau(X)$ is $C$-embedded in $M_\si(X)$, while $M_R(X)$ is not $C$-embedded in $M_\tau(X)$.

Our last example shows that the sufficient condition of Proposition~\ref{p:S} is not necessary for the $C$-embedding of $M_R(X)$ in $M(\be X)$. Recall that a set $A\subset[0,1]$ has universal measure zero if $\la^*(A)=0$ for every Borel probability measure $\la$ on $[0,1]$ vanishing on points. Uncountable sets of universal measure zero exist in ZFC: Grzegorek \cite{Grzegorek1980} constructed one of cardinality $\mathrm{non}(\mathcal N)$; see also \cite{Miller1984} for ZFC examples of cardinality $\aleph_1$. Every compact subset of a set $A$ of universal measure zero is countable, since an uncountable compact subset of $[0,1]$ contains a Cantor set, which carries a probability measure vanishing on points.

\begin{example}\label{ex:83}
Let $A\subset[0,1]$ be an uncountable set of universal measure zero, $T=[0,1]\times\{0,1\}^{\om_2}$, and for $c\in\{0,1\}$ let
\[
\Sigma_c=\set{y\in\{0,1\}^{\om_2}:|\set{\eta:y_\eta\ne c}|\le\om}.
\]
Let $X=(A\times\Sigma_1)\cup\big(([0,1]\setminus A)\times\Sigma_0\big)\subset T$, and let $f=\pi_0|_X\:X\to[0,1]$ be the projection to the first factor. Then
\begin{enumerate}
\item
$X$ is pseudocompact and $\be X=T$;
\item
$M_R(X)$ is $C$-embedded in $M(\be X)$;
\item
$f(X)=[0,1]$, but $f(S)\cap A$ is countable, and so $f(S)\ne f(X)$, for every $\si$-compact $S\subset X$.
\end{enumerate}
\end{example}
\begin{proof}
For $a\in[0,1]$ let $c(a)=1$ if $a\in A$ and $c(a)=0$ otherwise; thus $(a,y)\in X$ iff $y\in\Sigma_{c(a)}$. For $J\subset\om_2$ let $\pi_{0,J}\:T\to[0,1]\times\{0,1\}^J$ be the projection.

(i) Let $J\subset\om_2$ be countable, $a\in[0,1]$ and $u\in\{0,1\}^J$. The point $(a,y)$ with $y|_J=u$ and $y_\eta=c(a)$ for $\eta\notin J$ lies in $X$. Hence $\pi_{0,J}(X)=[0,1]\times\{0,1\}^J$ for every countable $J$, so $X$ is $G_\de$-dense in $T$ by Theorem~\ref{t:prod:main}, and (i) follows from Proposition~\ref{p:tau}.

(ii) We verify the condition of Theorem~\ref{t:dbR} for $p=R$ and $q=a$. Let $f'\:X\to K$ be a continuous map into a metrizable compact space, and let $\mu\in P(T)$. There are a countable $J\subset\om_2$ and a continuous $g$ with $\bar f'=g\circ\pi_{0,J}$ \cite[Problem 2.7.12(d)]{EngelkingBookRu}. Let $\rho=(\pi_0)_\sharp\mu\in P([0,1])$, and write $\rho=\sum_ir_i\de_{a_i}+\rho^c$, where the sum is countable and $\rho^c$ vanishes on points. Since $A$ has universal measure zero, there is a $\si$-compact $F\subset[0,1]\setminus A$ with $\rho^c([0,1]\setminus F)=0$. Let
\[
S=\big(F\times\{0,1\}^J\times\{\bar 0\}\big)\cup\bigcup_i\big(\{a_i\}\times\{0,1\}^J\times\{\bar c_i\}\big),
\]
where $\bar 0$ and $\bar c_i$ are the constant points of $\{0,1\}^{\om_2\setminus J}$ with values $0$ and $c(a_i)$, respectively. Then $S\subset X$ is $\si$-compact and $\pi_{0,J}(S)=(F\cup\{a_i\}_i)\times\{0,1\}^J$. Hence
\[
\bar f'^{-1}(f'(S))\supset\pi_{0,J}^{-1}(\pi_{0,J}(S))=\pi_0^{-1}(F\cup\{a_i\}_i),
\]
and $\mu(\bar f'^{-1}(f'(S)))\ge\rho(F\cup\{a_i\}_i)=1$. By Theorem~\ref{t:dbR}, $M_R(X)$ is $C$-embedded in $M(\be X)$.

(iii) Clearly $f(X)=[0,1]$. It suffices to show that $\pi_0(L)\cap A$ is countable for every compact $L\subset X$. Suppose not, and choose points $(a_\xi,y^\xi)\in L$, $\xi<\om_1$, with distinct $a_\xi\in A$. Then $y^\xi\in\Sigma_1$, so the sets $D_\xi=\set{\eta:y^\xi_\eta=0}$ are countable. Let $D=\bigcup_{\xi<\om_1}D_\xi$ and $E=\om_2\setminus D$; then $|E|=\aleph_2$. The set $L^1=\set{(a,y)\in L:y_\eta=1\text{ for all }\eta\in E}$ is compact and contains all the points $(a_\xi,y^\xi)$. If $(a,y)\in L^1$, then $y$ has $\aleph_2$ coordinates equal to $1$, so $y\notin\Sigma_0$, and hence $a\in A$. Thus $\pi_0(L^1)$ is an uncountable compact subset of $A$, a contradiction.
\end{proof}

\begin{note}\label{r:83}
In Example~\ref{ex:83} the $\si$-compact set $S$ in Theorem~\ref{t:dbR} \textup(for $p=R$, $q=a$\textup) must depend on $\mu$ even for the fixed map $f=\pi_0|_X$. Indeed, given a $\si$-compact $S\subset X$, take $a\in A\setminus\pi_0(S)$ and $\mu=\de_{(a,y)}$ with $(a,y)\in X$; then $\mu(\fb^{-1}(f(S)))=0$. Compare Remark~\ref{r:Sdep}, where $S$ depends on $f$. In \cite{RS2} it is shown that for this space $P_R(X)$ is even pseudocompact.
\end{note}

%======================================================================
\section{Questions}\label{sec:questions}

\begin{question}\label{q:1}
Find an internal characterization of the spaces $X$ for which $M_\tau(X)$ is $C$-embedded in $M_\si(X)$ \textup(equivalently, $P_\tau(X)$ is $G_\de$-dense in $P_\si(X)$\textup).
\end{question}

By Example~\ref{ex:tree}, realcompactness of $X$ is not sufficient, even for locally compact $X$; measure-compactness is not necessary, even for realcompact $X$ (Remark~\ref{r:Rc}). Sufficient conditions are given in Propositions~\ref{p:S} and~\ref{p:tau}. If $X$ is pseudocompact and locally compact, and all compact subsets of $X$ are countable \textup(for instance, $X$ is a pseudocompact Mr\'owka--Isbell space\textup), then $M_\tau(X)=M_d(X)$ and $M_\si(X)=M(\be X)$, so $M_\tau(X)$ is $C$-embedded in $M_\si(X)$ iff $\be X$ is scattered (Theorem~\ref{t:d}(ii)). Thus, already for pseudocompact spaces, the answer depends on the structure of $\be X$. A characterization in the spirit of Theorem~\ref{t:dbR} would require a lifting theorem for $\tau$-additive measures (Remark~\ref{r:tau}(ii)).

\begin{question}\label{q:groups}
Let $G$ be a pseudocompact topological group. Is $M_\tau(G)$ $C$-embedded in $M(\be G)$?
\end{question}

Recall that $\be G$ is a compact group in which $G$ is $G_\de$-dense \cite{ComfortRoss1966}. If $\be G$ is homeomorphic to a product of metrizable compact spaces, the answer is positive by Proposition~\ref{p:tau}. In the general case it is natural to replace the product measures in the proof of Proposition~\ref{p:tau} by the Haar measures of the cosets of closed normal subgroups $N$ of $\be G$ with $\be G/N$ metrizable.

%======================================================================


\begin{thebibliography}{99}

\bibitem{AngoaOrtiz-CastilloTamariz-Mascarua2014}
J.~Angoa, Y.~F.~Ortiz-Castillo, \'A.~Tamariz-Mascar\'ua,
\emph{Ultrafilters and properties related to compactness},
Topology Proc. \textbf{43} (2014), 183--200.

\bibitem{at2009}
A.~Arhangel'skii, M.~Tkachenko,
\emph{Topological Groups and Related Structures},
Atlantis Press, Paris, 2008.

\bibitem{BanakhChigogidzeFedorchuk2003}
T.~Banakh, A.~Chigogidze, V.~Fedorchuk,
\emph{On spaces of $\si$-additive probability measures},
Topology Appl. \textbf{133} (2003), no.~2, 139--155.

\bibitem{Bogachev2007Measure_theory}
V.~I.~Bogachev,
\emph{Measure Theory}, Vols.~I, II,
Springer, Berlin, 2007.

\bibitem{Bogachev2024-1ru}
V.~I.~Bogachev,
\emph{On compactification of spaces of measures},
Funct. Anal. Appl. \textbf{58} (2024), no.~1, 2--15; translated from Funktsional. Anal. i Prilozhen. \textbf{58} (2024), no.~1, 4--21.

\bibitem{Bogachev2024dan-ru}
V.~I.~Bogachev,
\emph{Compactification of spaces of measures and pseudocompactness},
Dokl. Math. \textbf{110} (2024), no.~1, 357--360; translated from Dokl. Ross. Akad. Nauk. Mat. Inform. Protsessy Upr. \textbf{518} (2024), 75--79.

\bibitem{ComfortRoss1966}
W.~W.~Comfort, K.~A.~Ross,
\emph{Pseudocompactness and uniform continuity in topological groups},
Pacific J. Math. \textbf{16} (1966), 483--496.

\bibitem{EngelkingBookRu}
R.~Engelking,
\emph{General Topology}, 2nd ed.,
Heldermann, Berlin, 1989.

\bibitem{Kakutani1943}
S.~Kakutani,
\emph{Notes on infinite product measure spaces. II},
Proc. Imp. Acad. Tokyo \textbf{19} (1943), 184--188.

\bibitem{KempermanMaharam1970}
J.~H.~B.~Kemperman, D.~Maharam,
\emph{$R^c$ is not almost Lindel\"of},
Proc. Amer. Math. Soc. \textbf{24} (1970), 772--773.

\bibitem{Knowles1967}
J.~D.~Knowles,
\emph{On the existence of non-atomic measures},
Mathematika \textbf{14} (1967), 62--67.

\bibitem{Korovin1992}
A.~V.~Korovin,
\emph{Continuous actions of pseudocompact groups and axioms of topological group},
Comment. Math. Univ. Carolin. \textbf{33} (1992), no.~2, 335--343.

\bibitem{GillmanJerison1960}
L.~Gillman, M.~Jerison,
\emph{Rings of Continuous Functions},
Van Nostrand, Princeton, 1960; reprinted as Grad. Texts in Math. \textbf{43}, Springer, New York, 1976.

\bibitem{Grzegorek1980}
E.~Grzegorek,
\emph{Solution of a problem of Banach on $\si$-fields without continuous measures},
Bull. Acad. Polon. Sci. S\'er. Sci. Math. \textbf{28} (1980), 7--10.

\bibitem{Hernandez-HernandezHrusak2018}
F.~Hern\'andez-Hern\'andez, M.~Hru\v s\'ak,
\emph{Topology of Mr\'owka--Isbell spaces},
in: Pseudocompact Topological Spaces, Springer, 2018, 253--289.

\bibitem{Miller1984}
A.~W.~Miller,
\emph{Special subsets of the real line},
in: Handbook of Set-Theoretic Topology, North-Holland, Amsterdam, 1984, 201--233.

\bibitem{Mrowka1977}
S.~Mr\'owka,
\emph{Some set-theoretic constructions in topology},
Fund. Math. \textbf{94} (1977), 83--92.

\bibitem{Mrowka1954}
S.~Mr\'owka,
\emph{On completely regular spaces},
Fund. Math. \textbf{41} (1954), 105--106.

\bibitem{PelczynskiSemadeni1959}
A.~Pe\l czy\'nski, Z.~Semadeni,
\emph{Spaces of continuous functions (III). Spaces $C(\Omega)$ for $\Omega$ without perfect subsets},
Studia Math. \textbf{18} (1959), 211--222.

\bibitem{RS}
E.~A.~Reznichenko, Yu.~V.~Sadovnichiy,
\emph{A new approach to universal measurability},
Vestnik Instituta Matematiki NAN KR (Bishkek), 2025, no.~2, 21--32.

\bibitem{RS0}
E.~A.~Reznichenko,
\emph{Dense subsets of locally convex spaces in the weak topology: $C$-embeddings and factorization of continuous functions},
submitted to Mat. Zametki.

\bibitem{RS2}
E.~A.~Reznichenko, Yu.~V.~Sadovnichiy,
\emph{$C^*$-embeddings of spaces of probability measures, nonnegative measures and unit balls},
in preparation.

\bibitem{Reznichenko2025betap}
E.~A.~Reznichenko,
\emph{Stone--\v Cech compactifications of spaces of probability measures},
Izv. Math. \textbf{90} (2026), no.~5, 993--1022; translated from Izv. Ross. Akad. Nauk Ser. Mat. \textbf{90} (2026), no.~5, 152--185.

\bibitem{Reznichenko1989}
E.~A.~Reznichenko,
\emph{A pseudocompact space in which only the subsets of full cardinality are not closed and not discrete},
Vestnik Moskov. Univ. Ser. I Mat. Mekh., 1989, no.~6, 69--70.

\bibitem{ReznichenkoTkachenko2024}
E.~A.~Reznichenko, M.~G.~Tkachenko,
\emph{All countable subsets of pseudocompact quasitopological Korovin groups are closed, discrete and $C^*$-embedded},
Topology Appl. \textbf{341} (2024), 108728.

\bibitem{ReznichenkoTkachenkoKorovin}
E.~A.~Reznichenko, M.~G.~Tkachenko,
\emph{Regular rigid Korovin orbits},
Rev. R. Acad. Cienc. Exactas F\'is. Nat. Ser. A Mat. RACSAM \textbf{120} (2026), no.~1, Paper No.~6.

\bibitem{Rudin1957}
W.~Rudin,
\emph{Continuous functions on compact spaces without perfect subsets},
Proc. Amer. Math. Soc. \textbf{8} (1957), 39--42.

\bibitem{Shakhmatov1986}
D.~Shakhmatov,
\emph{A pseudocompact Tychonoff space all countable subsets of which are closed and $C^*$-embedded},
Topology Appl. \textbf{22} (1986), no.~2, 139--144.

\bibitem{Wheeler1983}
R.~F.~Wheeler,
\emph{A survey of Baire measures and strict topologies},
Expo. Math. \textbf{1} (1983), 97--190.

\end{thebibliography}
\end{document}